\documentclass[12pt]{amsart}
\usepackage{graphicx}
\usepackage{amsmath}
\usepackage{amssymb} 
\usepackage[usenames,dvipsnames]{color}
\usepackage{pinlabel}

\newtheorem{thm}{Theorem}[section]  
\newtheorem{cor}[thm]{Corollary}

\newtheorem{defin}[thm]{Definition} 
\newtheorem{lemma}[thm]{Lemma} 
\newtheorem{example}[thm]{Example} 
\newtheorem{prop}[thm]{Proposition} 
\newtheorem{conj}{Conjecture} 
\newcommand{\aaa}{\mbox{$\alpha$}}
\newcommand{\bbb}{\mbox{$\beta$}}
\newcommand{\ddd}{\mbox{$\delta$}}

\newcommand{\bdd}{\mbox{$\partial$}}

\begin{document}  

\title{Fibered knots and Property 2R}   

\author{Martin Scharlemann}
\address{\hskip-\parindent
        Martin Scharlemann\\
        Mathematics Department\\
        University of California\\
        Santa Barbara, CA USA}
\email{mgscharl@math.ucsb.edu}

\author{Abigail Thompson}
\address{\hskip-\parindent
Abigail Thompson\\
Mathematics Department \\
University of California, Davis\\
Davis, CA 95616, USA}
\email{thompson@math.ucdavis.edu}

\thanks{Research partially supported by National Science Foundation grants.  Thanks also to Bob Gompf for helpful comments, and to Mike Freedman and Microsoft's Station Q for rapidly organizing a mini-conference on this general topic.}

\date{\today}

\begin{abstract} It is shown, using sutured manifold theory, that if there are any $2$-component counterexamples to the Generalized Property R Conjecture, then any knot of least genus among components of such counterexamples is not a fibered knot.  

The general question of what fibered knots might appear as a component of such a counterexample is further considered; much can be said about the monodromy of the fiber, particularly in the case in which the fiber is of genus two.
\end{abstract}

\maketitle

\section{Introductory remarks}  

Recall the famous Property R theorem, proven in a somewhat stronger form by David Gabai \cite{Ga2}:

\begin{thm}[Property R] \label{thm:PropR} If $0$-framed surgery on a knot $K \subset S^3$ yields $S^1 \times S^2$ then $K$ is the unknot. 
\end{thm}

There is a natural way of trying to generalize Theorem \ref{thm:PropR} to links in $S^3$.  In fact, there are several ways in which it can be generalized, but in the discussion here we restrict to the least complex one, known as the Generalized Property R Conjecture (see Conjecture \ref{conj:genR} below).  Other options are described in \cite{GSch} and briefly below.  The interest in this conjecture, as in the case of the original Property R Conjecture, is motivated in part by $4$-manifold questions. The viewpoint taken here derives almost entirely from 3-manifold theory, in particular new insights that can be brought to the question by developments in sutured manifold theory beyond those used by Gabai in his original proof of Property R.  

Unless explicitly stated otherwise, all manifolds throughout the paper will be compact and orientable.  

\section{Handle slides and Generalized Property R}  

To make sense of how Property R might be generalized, recall a small bit of $4$-manifold handlebody theory \cite{GS}.  Suppose $L$ is a link in a $3$-manifold $M$ and each component of $L$ is assigned a framing, that is a preferred choice of cross section to the normal bundle of the component in $M$. For example, if $M = S^3$, a framing on a knot is determined by a single integer, the algebraic intersection of the preferred cross-section with the longitude of the knot.  (In an arbitrary $3$-manifold $M$ a knot may not have a naturally defined longitude.)  Surgery on the link via the framing is standard Dehn surgery (though restricted to integral coefficients):  a regular neighborhood of each component is removed and then reattached so that the meridian is identified with the cross-section given by the framing.  Associated to this process is a certain $4$-manifold:  attach $4$-dimensional $2$-handles to $M \times I$ along $L \times \{ 1 \}$, using the given framing of the link components.  The result is a $4$-dimensional cobordism, called the {\em trace} of the surgery, between $M$ and the $3$-manifold $M'$ obtained by surgery on $L$. The collection of belt spheres of the $2$-handles constitute a link $L' \subset M'$ called the dual link; the trace of the surgery on $L \subset M$ can also be viewed as the trace of a surgery on $L' \subset M'$.

The $4$-manifold trace of the surgery on $L$ is unchanged if one $2$-handle is slid over another $2$-handle.  Such a handle slide is one of several moves allowed in the Kirby calculus \cite{Ki1}.  When the $2$-handle corresponding to the framed component $U$ of $L$ is slid over the framed component $V$ of $L$ the effect on the link is to replace $U$ by the band sum $\overline{U}$ of $U$ with a certain copy of $V$, namely the copy given by the preferred cross-section given by the framing of $V$.  

If $M$ is $S^3$ there is a simple formula for the induced framing on the new component $\overline{U}$.  Suppose $u, v \in \mathbb{Z}$ give the framing of the original components $U$ and $V$ respectively, and $U \cdot V \in \mathbb{Z}$ is the algebraic linking number of the components $U$ and $V$ in $S^3$.  Then the framing of the new component $\overline{U}$ that replaces $U$ is given by the formula  \cite[p.142]{GS}: $$u + v + 2 \; link(U, V).$$ Any statement about obtaining $3$-manifolds by surgery on a link will have to take account of this move, which we continue to call a handle-slide, in deference to its role in $4$-dimensional handle theory. 

Suppose $\overline{U} \subset M$ is obtained from components $U$ and $V$ by the handle-slide of $U$ over $V$ as described above.  Let $U' \subset M'$ and $V' \subset M'$ be the dual knots to $U$ and $V$.  It will be useful to note this counterintuitive but elementary lemma:

\begin{lemma}  \label{lemma:dual}  The link in $M'$ that is dual to $\overline{U} \cup V$ is $U' \cup \overline{V'}$, where $\overline{V'}$ is obtained by a handle-slide of $V'$ over $U'$.
\end{lemma}

\begin{proof}  It suffices to check this for the simple case in which the $3$-manifold is a genus $2$ handlebody, namely a regular neighborhood of $U$, $V$, and the arc between them along which the band-sum is done.  A sketch of this is shown in Figure \ref{fig:dual2}.  The dual knots $U' = \overline{U}'$, $V'$ and $\overline{V}'$ are displayed as boundaries of meridian disks for regular neighborhoods of $U$, $V$ and $\overline{V} = V$ respectively. 

 \begin{figure}[ht!]
 \labellist
\small\hair 2pt
\pinlabel \color{red}{$\overline{U}$} at 390 543
\pinlabel \color{red}{$U$} at 100 543
\pinlabel \color{black}{slide} at 171 512
\pinlabel \color{blue}{$\overline{V'}$} at 465 500
\pinlabel ${U'}$ at 90 435
\pinlabel ${\overline{U}' = U'}$ at 380 435
\pinlabel \color{ForestGreen}{$V$} at 245 545
\pinlabel \color{blue}{$V'$} at 225 435
 \color{black}
  \endlabellist
    \centering
    \includegraphics[scale=0.7]{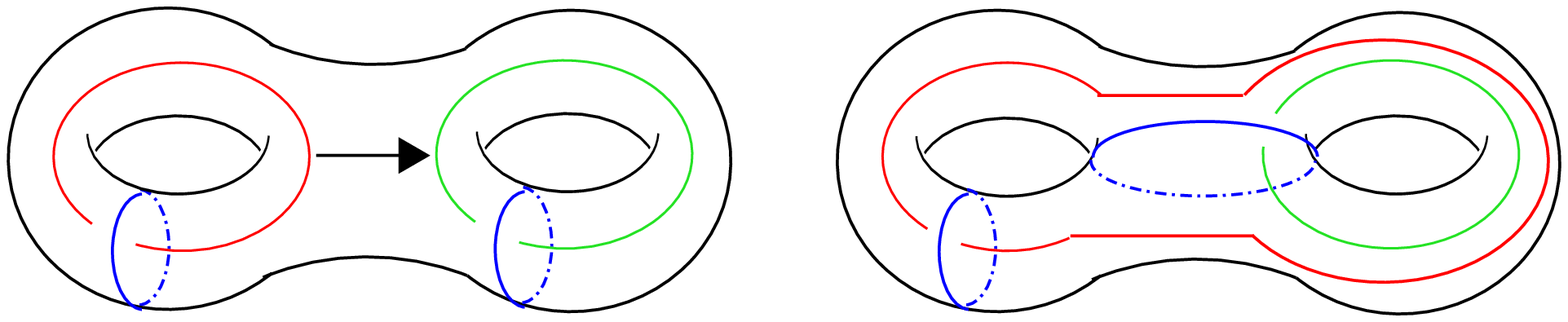}
    \caption{} \label{fig:dual2}
    \end{figure}

Alternatively, a $2$-dimensional schematic of the $4$-dimensional process is shown in Figure \ref{fig:dual}.  The handle corresponding to $U$ is shown half-way slid across the handle corresponding to $V$.  Each disk in the figure is given the same label as its boundary knot in $M$ or $M'$ as appropriate. \end{proof}

     \begin{figure}[ht!]
 \labellist
\small\hair 2pt
\pinlabel \color{red}{$\overline{U}$} at 37 183
\pinlabel \color{red}{$U$} at 137 86
\pinlabel \color{black}{slide} at 346 168
\pinlabel \color{blue}{$\overline{V'}$} at 238 122
\pinlabel ${U' = \overline{U}'}$ at 162 238
\pinlabel \color{ForestGreen}{$V$} at 310 50
\pinlabel \color{blue}{$V'$} at 396 101
 \color{black}
  \endlabellist
    \centering
    \includegraphics[scale=0.6]{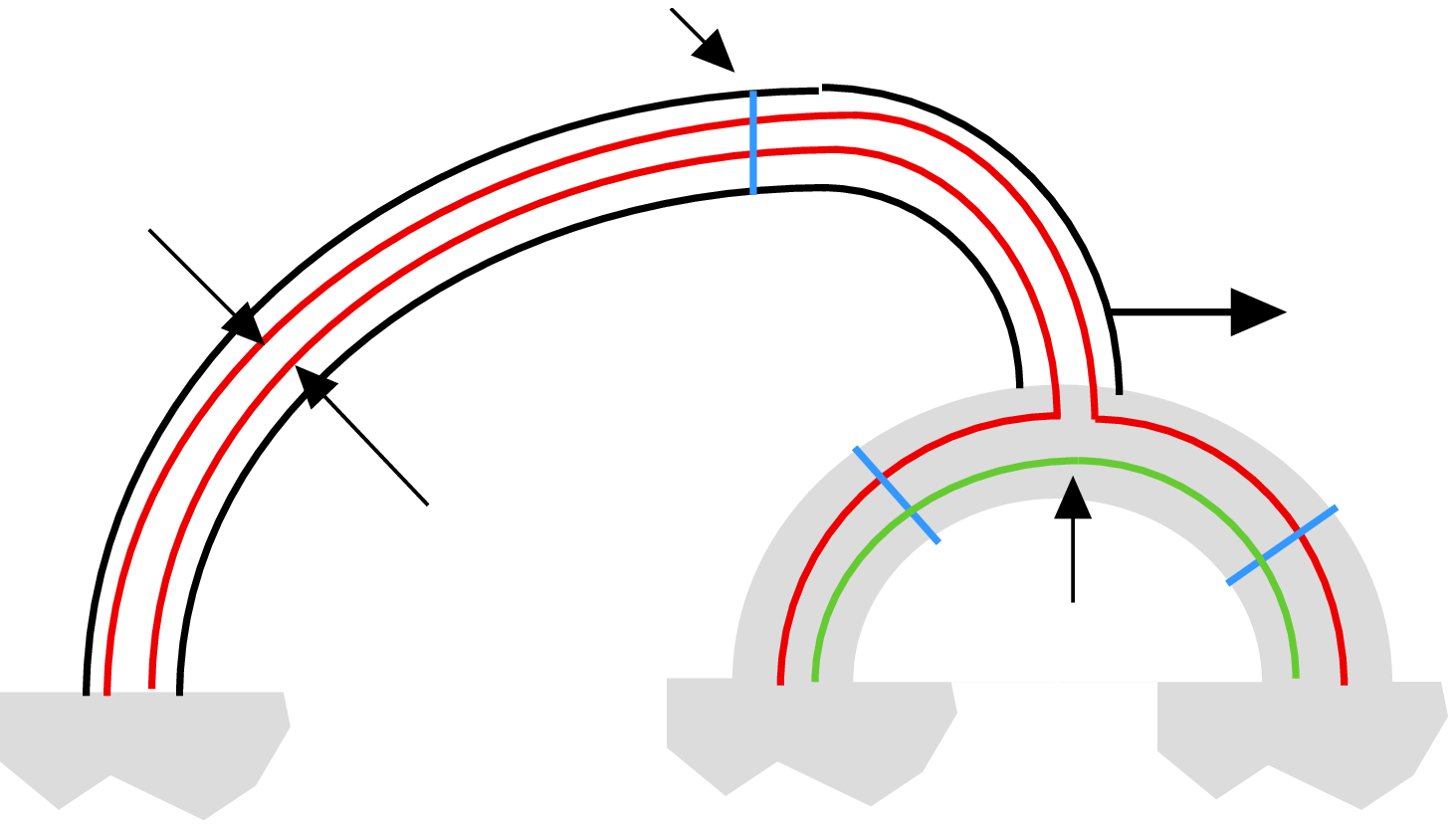}
    \caption{} \label{fig:dual}
    \end{figure}

Let $\#_{n} (S^{1} \times S^{2})$ denote the connected sum of $n$ copies of $S^1 \times S^2$.  The Generalized Property R conjecture (see \cite[Problem 1.82]{Ki2}) says this:  

\begin{conj}[Generalized Property R] \label{conj:genR} Suppose $L$ is an integrally framed link of $n \geq 1$ components in $S^3$, and surgery on $L$ via the specified framing yields $\#_{n} (S^{1} \times S^{2})$.  Then there is a sequence of handle slides on $L$ that converts $L$ into a $0$-framed unlink.  
\end{conj}

\begin{figure}[ht!]
 \labellist
\small\hair 2pt
\pinlabel $0$ at 72 89
\pinlabel $0$ at 137 120
\pinlabel $0$ at 281 18
\pinlabel $0$ at 338 120
\pinlabel $0$ at 432 18
\pinlabel $0$ at 432 120
\pinlabel {band sum here} at 108 -11

 \endlabellist
    \centering
    \includegraphics[scale=0.7]{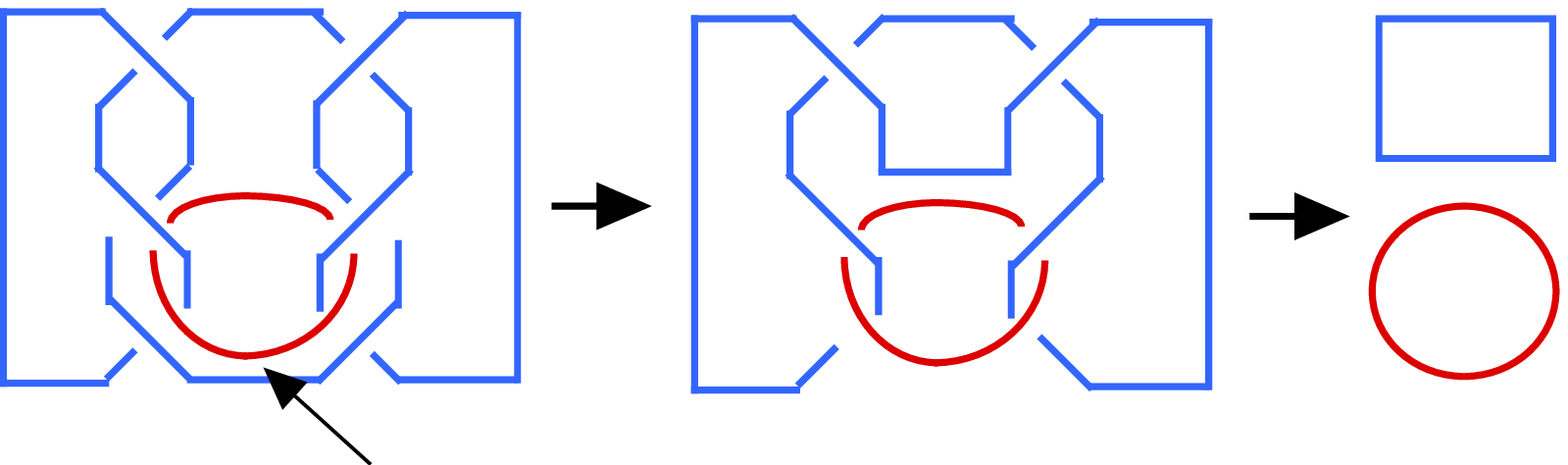}
    \caption{} \label{fig:squareknot}
    \end{figure}

\begin{figure}[ht!]
    \centering
    \includegraphics[scale=0.7]{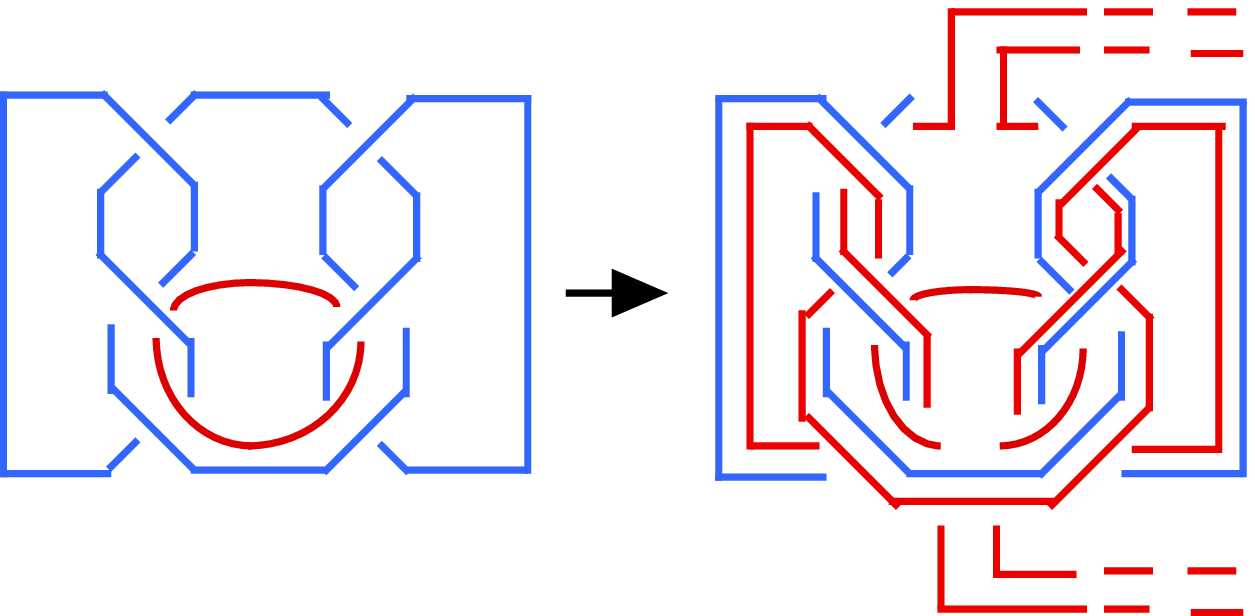}
    \caption{} \label{fig:squareknot2}
    \end{figure}

In the case $n = 1$ no slides are possible, so Conjecture \ref{conj:genR} does indeed directly generalize Theorem \ref{thm:PropR}.  On the other hand, for $n > 1$ it is certainly necessary to include the possibility of handle slides.  Figure \ref{fig:squareknot} shows an example of a more complicated link  on which $0$-framed surgery creates $\#_{2} (S^{1} \times S^{2})$.  To see this, note that the Kirby move shown, band-summing the square knot component to a copy of the unknotted component,  changes the original link to the unlink of two components, for which we know surgery yields $\#_{2} (S^{1} \times S^{2})$.  Even more complicated links with this property can be obtained, simply by using Kirby moves that complicate the link rather than simplify it.  See Figure \ref{fig:squareknot2}; the free ends of the band shown can be connected in an arbitrarily linked or knotted way.

The conjecture can be clarified somewhat by observing that the only framed links that are relevant are those in which all framings and linking numbers are trivial.  There is a straightforward $4$-dimensional proof, using the intersection pairing on the trace of the surgery.  Here is an equally elementary $3$-dimensional proof:

\begin{prop} \label{prop:frame}  Suppose $L$ is a framed link of $n \geq 1$ components in $S^3$, and surgery on $L$ via the specified framing yields $\#_{n} (S^{1} \times S^{2})$.  Then the components of $L$ are algebraically unlinked and the framing on each component is the $0$-framing. 
\end{prop}  

\begin{proof}  It follows immediately from Alexander duality that $H_1(S^3 - \eta(L)) \cong H^1(\eta(L)) \cong n\mathbb{Z}.$  In particular, filling in the solid tori via whatever framing we are given yields an epimorphism, hence an isomorphism $H_1(S^3 - \eta(L)) \to H_1(\#_{n} (S^{1} \times S^{2}))$.  For each torus component $T$ of $\bdd \eta(L)$, the filling will kill some generator of $H_1(T)$, so the homomorphism $H_1(T) \to H_1(\#_{n} (S^{1} \times S^{2}))$ is not injective. It follows that the homomorphism $H_1(T) \to H_1(S^3 - \eta(L))$ cannot be injective and, moreover, $ker(H_1(T) \to H_1(S^3 - \eta(L)))$ must contain the framing curve.  But $ker(H_1(T) \to H_1(S^3 - \eta(L)))$ must be contained in the subgroup generated by the standard longitude, since this is the only subgroup that is trivial when we just replace all the other components of $\eta(L)$.  It follows that the framing at each component is that of the standard longitude, namely the $0$-framing.   Since the longitude of each $T$ is null-homologous in $H_1(S^3 - \eta(L))$ it follows that all linking numbers are trivial.  
\end{proof} 

There is also an immediate topological restriction on the link itself, which carries over to a restriction on the knots can appear as individual components of such a link:

\begin{prop} \label{prop:slice}  Suppose $L$ is a framed link of $n \geq 1$ components in $S^3$, and surgery on $L$ via the specified framing yields $\#_{n} (S^{1} \times S^{2})$.  Then $L$ bounds a collection of $n$ smooth $2$- disks in a $4$-dimensional homotopy ball bounded by $S^3$. 
\end{prop}  

An equivalent way of stating the conclusion, following Freedman's proof of the $4$-dimensional topological Poincare Conjecture \cite{Fr}, is that $L$ (and so each component of $L$) is topologically slice in $B^4$.  

\begin{proof}  Consider the $4$-manifold trace $W$ of the surgery on $L$.   $\bdd W$ has one end diffeomorphic to $S^3$ and the other end, call it $\bdd_1 W$, diffeomorphic to $\#_{n} (S^{1} \times S^{2})$.  $W$ has the homotopy type of a once-punctured $\natural_{n} (B^{2} \times S^{2})$.  Attach $\natural_{n} (S^{1} \times B^{3})$ to $\bdd_1 W$ via the natural identification $\bdd B^{3} \cong S^2$.  The result is a homotopy $4$-ball, and the cores of the original $n$ $2$-handles that are attached to $L$ are the required $n$ $2$-disks.
\end{proof}

Somewhat playfully, we can turn the Generalized Property R Conjecture, which is a conjecture about links, into a conjecture about knots, and also stratify it by the number of components, via the following definition and conjecture.  

\begin{defin} A knot $K \subset S^3$ has {\bf Property nR} if it does not appear among the components of any $n$-component counterexamples to the Generalized Property R conjecture.
\end{defin}

\begin{conj}[Property nR Conjecture] All knots have Property nR.
\end{conj}

Thus the Generalized Property R conjecture for all $n$ component links is equivalent to the Property nR Conjecture for all knots.  Following Proposition \ref{prop:slice} any non-slice knot has Property nR for all $n$.  The first thing that we will show (Theorem \ref{thm:main}) is that if there are any counterexamples to Property 2R, a least genus such example cannot be fibered.  We already know that both of the genus one fibered knots (the figure 8 knot and the trefoil) cannot be counterexamples, since they are not slice. So these simplest of fibered knots do have Property 2R.  On the other hand, for reasons connected to difficulties proving the Andrews-Curtis Conjecture, there is strong evidence (see \cite{GSch}) that Property 2R may fail for as simple a knot as the square knot.  Since the square knot is fibered, it would then follow from Theorem \ref{thm:main} that there is a counterexample to Property 2R among genus one knots.  

\section{Special results for Property 2R}

Almost nothing is known about Generalized Property R, beyond the elementary facts noted in Propositions \ref{prop:frame} and \ref{prop:slice} that the framing and linking of the components of the link are all trivial and the links themselves are topologically slice.
A bit more is known about Property 2R.  The first was shown to us by Alan Reid:

\begin{prop}[A. Reid]  Suppose $L \subset S^3$ is a $2$-component link with tunnel number $1$.  If surgery on $L$ gives $\#_{2} (S^{1} \times S^{2})$ then $L$ is the unlink of two components.
\end{prop}

\begin{proof} The assumption that $L$ is tunnel number $1$ means that there is a properly embedded arc $\aaa \subset S^3 - \eta(L)$ so that  $S^3 - (\eta(L) \cup \eta(\aaa))$ is a genus $2$ handlebody $H$.  Let $G = \pi_1(S^3 - \eta(L))$.  There is an obvious epimorphism $Z*Z \cong \pi_1(H) \to G$ (fill in a meridian disk of $\eta(\aaa)$) and an obvious epimorphism $G \to \pi_1(\#_{2} (S^{1} \times S^{2})) \cong Z*Z$ (fill in solid tori via the given framing).  But any epimorphism $Z*Z \to Z*Z$ is an isomorphism, since free groups are Hopfian, so in fact $G \cong Z*Z$.  It is then a classical result that $L$ must be the unlink on two components. \end{proof} 

This first step towards the Property 2R conjecture is a bit disappointing, however, since handle-slides (the new and necessary ingredient for Generalized Property R) do not arise.  In contrast, Figure \ref{fig:squareknot} shows that handle slides are needed in the proof of the following:

\begin{prop}  \label{prop:unknot}  The unknot has Property 2R.
\end{prop}

\begin{proof}  Suppose $L$ is the union of two components,  the unknot $U$ and another knot $V$, and suppose some surgery on $L$ gives $\#_{2} (S^{1} \times S^{2})$.  Following Proposition \ref{prop:frame} the surgery is via the $0$-framing on each and, since $U$ is the unknot, $0$-framed surgery on $U$ alone creates $S^1 \times S^2$.  Moreover, the curve $U' \subset S^1 \times S^2$ that is dual to $U$ is simply $S^1 \times \{ p \}$ for some point $p \in S^2$.  

A first possibility is that $V$ is a satellite knot in $S^1 \times S^2$, so $V$ lies in a solid torus $K$ in such a way that the torus $\bdd K$ is essential in $S^1 \times S^2 - \eta(V)$.  Since there is no essential torus in $\#_{2} (S^{1} \times S^{2})$, $\bdd K$ compresses after the surgery on $V$.  Since $\#_{2} (S^{1} \times S^{2})$ contains no summand with finite non-trivial homology, it follows from the main theorem of \cite{Ga3} that $V$ is a braid in $K$ and that surgery on $V$ has the same effect on $S^1 \times S^2$ as some surgery on $K$.  Proposition \ref{prop:frame} shows that the surgery on $K$ must be along a longitude of $K$, but that would imply that $V$ has winding number $1$ in $K$.  The only braid in a solid torus with winding number $1$ is the core of the solid torus, so in fact $V$ is merely a core of $K$ and no satellite.  So we conclude that $V \subset S^1 \times S^2$ cannot be a satellite knot.

Consider the manifold $M = S^1 \times S^2 - \eta(V)$.  If $M$ is irreducible, then it is a taut sutured manifold (see, for example, \cite{Ga1}) and two different fillings (trivial vs. $0$-framed) along $\bdd \eta(V)$ yield reducible, hence non-taut sutured manifolds.  This contradicts \cite{Ga1}.  We conclude that $M$ is reducible.  It follows that $V$ is isotopic \underline{in $M$} to a knot $\overline{V}$ lying in a $3$-ball in $M - U'$ and that surgery on $\overline{V} \subset B^3$ creates a summand of the form $S^1 \times S^2$.  By Property R, we know that $\overline{V}$ is the unknot in $B^3$.  Hence $U \cup \overline{V} \subset S^3$ is the unlink of two components.

The proof, though, is not yet complete, because the isotopy of $V$ to $\overline{V}$ in $M$ may pass through $U'$.  But passing $V$ through $U' \subset M$ can be viewed as band-summing $V$ to the boundary of a meridian disk of $U'$ in $M$.  So the effect in $S^3$ is to replace $V$ with the band sum of $V$ with a longitude of $U$.  In other words, the knot $\overline{V}$, when viewed back in $S^3$, is obtained by from $V$ by a series of handle slides over $U$, a move that is permitted under Generalized Property R.  
\end{proof}

In a similar spirit, the first goal of the present paper is to prove a modest generalization of Proposition \ref{prop:unknot}.   A pleasant feature is that, since the square knot is fibered, Figure \ref{fig:squareknot2} shows that the proof will require handle slides of {\em both} components of the link.  

\begin{thm} \label{thm:main} No smallest genus counterexample to Property 2R is fibered.
\end{thm}

\begin{proof}  Echoing the notation of Proposition \ref{prop:unknot}, suppose there is a $2$-component counterexample to Generalized Property R consisting of a fibered knot $U$ and another knot $V$.  Let $M$ be the $3$-manifold obtained by $0$-framed surgery on $U$ alone.  Since $U$ is a fibered knot, $M$ fibers over the circle with fiber $F$, a closed orientable surface of the same genus as $U$.  The dual to $U$ in $M$ is a knot $U'$ that passes through each fiber exactly once.  

The hypothesis is that $0$-framed surgery on $V \subset M$ creates $\#_{2} (S^{1} \times S^{2})$.  Following  \cite[Corollary 4.2]{ST}, either  the knot $V$ lies in a ball, or $V$ is cabled with the surgery slope that of the cabling annulus, or $V$ can be isotoped in $M$ to lie in a fiber, with surgery slope that of the fiber.  If $V$ were cabled, then the surgery on $K$ would create a Lens space summand, which is clearly impossible in $\#_{n} (S^{1} \times S^{2})$.  If $V$ can be isotoped into a ball or into a fiber, then, as argued in the proof of Proposition \ref{prop:unknot}, the isotopy in $M$ is realized in $S^3$ by handle-slides of $V$ over $U$, so we may as well regard $V$ as lying either in a ball that is disjoint from $U'$ or in a fiber $F_0 \subset M$.  The former case, $V$ in a ball disjoint from $U'$ would, as in Proposition \ref{prop:unknot}, imply that the link $U \cup V \subset S^3$ is the unlink.  So we can assume that $V \subset F_0 \subset M$.

The surgery on $V$ that changes $M$ to $\#_{2} (S^{1} \times S^{2})$ has this local effect near $F_0$:  $M$ is cut open along $F_0$, creating two copies $F^{\pm}_0$, a $2$-handle is attached to  the copy of $V$ in each of $F^{\pm}_0$, compressing the copies of the fiber to surfaces $F'^{\pm}$.  The surfaces $F'^{\pm}$ are then glued back together by the obvious identification to give a surface $F' \subset \#_{2} (S^{1} \times S^{2})$. (See the Surgery Principle Lemma \ref{lemma:surgprin} below for more detail.)   This surface has two important features:  each component of $F'$ (there are two components if and only if $V$ is separating in $F$) has lower genus than $F$; and $F'$ intersects $U'$ in a single point.  

Let $V' \subset \#_{2} (S^{1} \times S^{2})$ be the dual knot to $V$ and let $F''$ be the component of $F'$ that intersects $U'$.  $V'$ intersects $F'$ in some collection of points (in fact, two points, but that is not important for the argument).  Each point in $V' \cap F''$ can be removed by a handle-slide of $V'$ over $U'$ along an arc in $F''$.  Let $V''$ be the final result of these handle-slides.  Then $F''$ is an orientable surface that has lower genus than $F$, is disjoint from $V''$ and intersects $U'$ in a single point.

Following Lemma \ref{lemma:dual} the handle-slides of $V'$ over $U'$ in $\#_{2} (S^{1} \times S^{2})$ correspond in $S^3$ to handle-slides of $U$ over $V$.  Call the knot in $S^3$ that results from all these handle-slides $\overline{U} \subset S^3$.  Since $F''$ is disjoint from $V''$, and intersects $U'$ in a single meridian, $F'' - U'$ is a surface in $S^3 - \overline{U}$ whose boundary is a longitude of $\overline{U}$.  In other words, the knot $\overline{U}$, still a counterexample to Property 2R, has $$genus(\overline{U}) = genus(F'') < genus(F) = genus(U)$$ as required.
\end{proof} 

\section{Fibered manifolds and Heegaard splittings}

We have just seen that a fibered counterexample to Property 2R would not be a least genus counterexample.  We now explore other properties of potential fibered counterexamples.  In this section we consider what can be said about the monodromy of a fibered knot in $S^3$, and the placement of a second component with respect to the fibering, so that surgery on the $2$-component link yields $\#_{2} (S^{1} \times S^{2})$.  Perhaps surprisingly, the theory of Heegaard splittings is useful in answering these questions.  Much of this section in fact considers the more general question of when $\#_{2} (S^{1} \times S^{2})$ can be created by surgery on a knot in a $3$-manifold $M$ that fibers over a circle.  The application to Property 2R comes from the special case in which the manifold $M$ is obtained from $0$-framed surgery on a fibered knot in $S^3$.  

Suppose $F$ is a surface in a $3$-manifold $M$ and $c \subset F$ is an essential simple closed curve in $F$.  A tubular neighborhood $\eta(c) \subset M$ intersects $F$ in an annulus; the boundary of the annulus in $\bdd \eta(c)$ defines a slope on $\eta(c)$.  Let $M_{surg}$ denote the manifold obtained from $M$ by surgery on $c$ with this slope and let $F'$ be the surface obtained from $F$ by compressing $F$ along $c$.   

\begin{lemma}[Surgery Principle] \label{lemma:surgprin}  $M_{surg}$ can be obtained from $M$ by the following $3$-step process:
\begin{enumerate}
\item Cut $M$ open along $F$, creating two new surfaces $F_{\pm}$ in the boundary, each homeomorphic to $F$.   
\item Attach a $2$-handle to each of $F_{\pm}$ along the copy of $c$ it contains.  This changes each of the new boundary surfaces $F_{\pm}$ to a copy of $F'$.  Denote these two surfaces $F'_{\pm}$.
\item Glue $F'_+$ to $F'_-$ via the natural identification.
\end{enumerate}
\end{lemma}

\begin{proof}  The surgery itself is a $2$-step process: Remove a neighborhood of $\eta(c)$, then glue back a copy of $S^1 \times D^2$ so that $\{point \} \times \bdd D^2$ is the given slope.  The first step is equivalent to cutting $F$ along an annulus neighborhood $A$ of $c$ in $F$, creating a torus boundary component as the boundary union of the two copies $A_{\pm}$ of $A$.  Thinking of $S^1$ as the boundary union of two intervals, the second step can itself be viewed as a two-step process:  attach a copy of $I \times D^2$ to each annulus $A_{\pm}$ along $I \times \bdd D^2$ (call the attached copies $(I \times D^2)_{\pm}$), then identify the boundary disks $(\bdd I \times D^2)_+$ with $(\bdd I \times D^2)_-$ in the natural way.  This creates a three-stage process which is exactly that defined in the lemma, except that in the lemma $F-A$ is first cut apart and then reglued by the identity. \end{proof}

The case in which $M$ fibers over a circle with fiber $F$ is particularly relevant.  We will stay in that case throughout the remainder of this section (as always, restricting to the case that $M$ and $F$ are orientable) and use the following notation:
\begin{enumerate} 
\item $h: F \to F$ is the monodromy homeomorphism of $M$.
\item $c$ is an essential simple closed curve in $F$.
\item $F'$ is the surface obtained by compressing $F$ along $c$
\item $M_{surg}$ is the manifold obtained by doing surgery on $M$ along $c \subset F \subset M$ using the framing given by $F$.
\end{enumerate}
Note that $F'$ may be disconnected, even if $F$ is connected.

\begin{prop}  \label{prop:isotopic1}  Suppose $h(c)$ is isotopic to $c$ in $F$
\begin{itemize}
\item If $c$ is non-separating in $F$, or if $c$ is separating and the isotopy from $h(c)$ to $c$ reverses orientation of $c$, then $M_{surg} \cong N \# (S^1 \times S^2)$, where $N$ fibers over the circle with fiber $F'$.
\item If $c$ separates $F$ so $F' = F_1 \cup F_2$, and the isotopy from $h(c)$ to $c$ preserves orientation of $c$, then $M_{surg} \cong M_1 \# M_2$, where each $M_i$ fibers over the circle with fiber $F_i$.
\end{itemize}
\end{prop}

\begin{proof}  We may as well assume that $h(c) = c$ and consider first the case where $h|c$ is orientation preserving.  In this case, the mapping cylinder of $c$ in $M$ is a torus $T$ containing $c$.  The $3$-stage process of Lemma \ref{lemma:surgprin} then becomes:  
\begin{enumerate}
\item $M$ is cut along $T$ to give a manifold $M_-$ with two torus boundary components.  $M_-$ fibers over the circle with fiber a twice-punctured $F'$. ($F'$ is connected if and only if $c$ is non-separating.)  
\item A $2$-handle is attached to each torus boundary component $T_{\pm}$, turning the boundary into two $2$-spheres.  
\item The two $2$-spheres are identified.
\end{enumerate}
The second and third stage together are equivalent to filling in a solid torus along each $T_{\pm}$, giving an $F'$-fibered manifold $M'$, then removing a $3$-ball from each solid torus and identifying the resulting $2$-spheres.  Depending on whether $F'$ is connected or not, this is equivalent to either adding $S^1 \times S^2$ to $M'$ or adding the two components of $M'$ together.

The case in which $h|c$ is orientation reversing is only slightly more complicated.  Since $M$ is orientable, the mapping cylinder of $h|c$ is a $1$-sided Klein bottle $K$, so $\bdd(M - \eta(K))$ is a single torus $T$.  The argument of Lemma \ref{lemma:surgprin} still mostly applies, since $c$ has an annulus neighborhood in $K$, and shows that the surgery can be viewed as attaching two $2$-handles to $T$ along parallel curves, converting the boundary into two $2$-spheres, then identifying the $2$-spheres.  This is again equivalent to filling in a solid torus at $T$ (which double-covers $S^1$) and then adding $S^1 \times S^2$.  But filling in a solid torus at $T \subset (M - \eta(K))$ changes the fiber from $F$ to $F'$.  (Note that if $c$ separates $F$, so $F' = F_1 \cup F_2$, then since $h$ is orientation preserving on $F$ but orientation reversing on $c$, $h$ must exchange the $F_i$.  So $N$ also fibers over the circle with fiber $F_1$.)
\end{proof}

\begin{cor}  \label{cor:isotopic2} If $M_{surg} \cong \#_{2} (S^{1} \times S^{2})$ and $h(c)$ is isotopic to $c$ in $F$, then $F$ is a torus.
\end{cor}

\begin{proof}   According to Proposition \ref{prop:isotopic1}, the hypotheses imply that $S^1 \times S^2$ fibers over the circle with fiber (a component of) $F'$.  But this forces $F' \cong S^2$ and so $F \cong T^2$. \end{proof}

Surgery on fibered manifolds also provides a natural connection between the surgery principle and Heegaard theory:

\begin{defin}  Suppose $H_1$, $H_2$ are two copies of a compression body $H$ and $h: \bdd_+ H  \to \bdd_+ H$ is a given homeomorphism.  Then the union of $H_1$, $H_2$ along their boundaries, via the identity on $\bdd_- H_i$ and via $h:  \bdd_+ H_1  \to \bdd_+ H_2$, is called the {\em Heegaard double} of $H$ (via $h$).
\end{defin}

Lemma \ref{lemma:surgprin} gives this important example:

\begin{example}  \label{example:double}  For $M, F, h, c, M_{surg}$ as above, let $H$ be the compression body obtained by attaching a $2$-handle to $F \times \{ 1 \} \subset F \times I$ along $c$.  Then $M_{surg}$ is the Heegaard double of $H$ via $h$.
\end{example}

Note that the closed complement $N$ of $\bdd_- H_1 = \bdd_- H_2$ in any Heegaard double is a manifold with Heegaard splitting $N \cong H_1 \cup_{\bdd_+} H_2$.  Here is a sample application, using Heegaard theory:

\begin{prop}  \label{prop:monodromy} For $M, F, h, c, M_{surg}$ as above, suppose some surgery on $c$ gives a reducible manifold.  Then the surgery slope is that of $F$ and either 
\begin{enumerate}
\item $h(c)$ can be isotoped in $F$ so that it is disjoint from $c$ or 
\item $c \subset F$ is non-separating and $M_{surg} \cong N \# L$, where 
\begin{itemize}
\item $N$ fibers over the circle with fiber $F'$ and 
\item $L$ is either $S^3$ or a Lens space.
\end{itemize}
\end{enumerate}

\end{prop}

Note in particular that possibility (2) is not consistent with $M_{surg} \cong \#_{2} (S^{1} \times S^{2})$. 

\begin{proof}   Choose distinct fibers $F_0, F_1$ in $M$, with $c \subset F_1$.  Via \cite[Corollary 4.2]{ST} and the proof of Theorem \ref{thm:main} we know that the surgery on $c$ must use the framing given by the fiber $F_1$, so the result of surgery is $M_{surg}$.   Example \ref{example:double} shows that $M_{surg}$ is a Heegaard double via $h$, so the complement $M_- = M_{surg} - \eta(F')$ of a regular neighborhood of $F' = \bdd_- H$ has a Heegaard splitting $H_1 \cup_{F_0} H_2$.  That is,  $F_0 = \bdd_+ H_1 = \bdd_+ H_2$.

If $F' \cong S^2$, so $F \cong T^2$, then $M_{surg} \cong M_- \# S^1 \times S^2$.  Since $F \cong T^2$, the Heegaard splitting $H_1 \cup_{F_0} H_2$ of $M_-$ is of genus $1$, so $M_-$ is either $S^3$, a Lens space, or $S^1 \times S^2$.  But the last happens only if the same curve in $F_0$ compresses in both $H_1$ and $H_2$; in our context, that implies $c$ and $h(c)$ are isotopic in $F$, and so can be isotoped to be disjoint.

If $F' \ncong S^2$, choose a reducing sphere with a minimal number of intersection curves with $F'$.  If the reducing sphere is disjoint from $F'$, then $M_-$ is reducible. If the reducing sphere intersects $F'$, then at least one copy of $F'$ in $\bdd M_-$ must be compressible in $N$.  We conclude that in either case the Heegaard splitting $H_1 \cup_{F_0} H_2$ of $M_-$ is weakly reducible (and possibly reducible), see \cite{CG}.  That is, there are essential disjoint simple closed curves $\aaa_1, \aaa_2$ in $F = \bdd_+ H_i$ which compress respectively in $H_1$ and $H_2$. 
 
 \medskip 

{\bf Case 1:}  The curve $c$ is separating.

In this case, since the compression bodies $H_i$ each have only the $2$-handle with boundary $c \subset F_1$ attached, any curve in $\bdd _+ H_i = F_0$ that compresses in $H_i$ is isotopic to $c \subset \bdd_+ H_i \cong F$.  In particular, fixing the identification $F_0 = \bdd_+ H_2$, $\aaa_2$  must represent $c$ in $F_0$ and $\aaa_1$ represents $h(c)$.  Hence $c$ and $h(c)$ are disjoint.

\medskip

{\bf Case 2:}  The curve $c$ is non-separating, and so is at least one of the curves $\aaa_1, \aaa_2$.

If both curves $\aaa_i$ are non-separating then, as in Case 1, $\aaa_1$ and $\aaa_2$, when viewed in the handlebodies $H_1, H_2$, must each be isotopic to $c \subset \bdd_+ H_i \cong F_0$ and the case concludes as Case 1 did.

If $\aaa_2$ is non-separating, and $\aaa_1$ is separating, then $\aaa_2$  is isotopic to $c \subset \bdd_+ H_2 = F_0$ whereas $\aaa_1$ bounds a punctured torus $T \bdd_+ H_2$ on which $h(c)$ lies.  If $\aaa_2$ is disjoint from $T$, then $c$ and $h(c)$ are disjoint, as required.  If $\aaa_2$ lies in $T$ then $\bdd T$ also bounds a disk in $H_2$.  The union of the disks in $H_1$ and $H_2$ bounded by $\bdd T$ is a sphere that decomposes $M_-$ into $F' \times I \# W$.  This implies that $M_{surg} \cong N \# W$, where $N$ fibers over $S^1$ with fiber $F'$ and $W$ is Heegaard split by $T$ into two solid tori, with meridian disks bounded by $c$ and $h(c)$ respectively.  If $|c \cap h(c)| > 1$ then $W$ is a Lens space. If $|c \cap h(c)| = 1$ then $W = S^3$.  If $|c \cap h(c)| = 0$ then $h(c)$ is disjoint from $c$.

\medskip

{\bf Case 3:} The curve $c$ is non-separating, but both $\aaa_1, \aaa_2$ are separating.

In this case, much as in Case 2, each $\aaa_i$ cuts off a torus $T_i$ from $\bdd_+ H_2 = F_0$, with $c \subset T_2$ and $h(c) \subset T_1$.  Since the $\aaa_i$ are disjoint, the two tori either also are disjoint (and the proof is complete) or the two tori coincide.  If the two tori coincide, the argument concludes as in Case 2.  \end{proof}

\section{Could there be fibered counterexamples of genus $2$?}

In applying Proposition \ref{prop:monodromy} to the manifold $M$ obtained from $0$-framed surgery on a fibered knot $K \subset S^3$, note that the isotopy in the Proposition takes place in a fiber $F$ of $M$, the closed manifold obtained by $0$-framed surgery on $K$, not in the fiber $F - \{ point \}$ of the knot $K$ itself.   The importance of the distinction is illustrated by the following Proposition which, without the distinction, would (following Propositions \ref{prop:unknot} and \ref{prop:monodromy}) seem to guarantee that all genus $2$ fibered knots have Property 2R.

\begin{prop}  Suppose $U \subset S^3$ is a fibered knot, with fiber the punctured surface $F_- \subset S^3$ and monodromy $h_-: F_- \to F_-$.  Suppose a knot $V \subset F_-$ has the property that $0$-framed surgery on the link $U \cup V$ gives $\#_{2} (S^{1} \times S^{2})$ and $h_-(V)$ can be isotoped to be disjoint from $V$ in $F_-$.  Then either $V$ is the unknot or $genus(F_-) \neq 1, 2$.  \end{prop}

\begin{proof}  {\em Case 1:} $V$ bounds a disk in $F_-$ or is parallel in $F_-$ to $\bdd F_- = U$.

In this case, $0$-framed surgery on $U \cup V$ would be $N \# S^1 \times S^2$, where $N$ is the result of $0$-framed surgery on $U$.  Our hypothesis is that $N \cong S^1 \times S^2$ which, by classical Property R \cite{Ga2}, implies that $U$ is the unknot.  Hence $genus(F_-) = 0$.  

\bigskip

{\em Case 2:} $V$ is essential in $F_-$.  

If $F_-$ is a punctured torus, then the fact that $V$ is essential and $h_-(V)$ can be isotoped off of $V$ imply that $h_-(V)$ is isotopic to $V$, and we may as well assume that $h_-(V) = V$.  The mapping torus of $h_-|V$ is then a non-separating torus in $S^3$, which is absurd.  

Suppose $F_-$ is a punctured genus $2$-surface, and let $F$ denote the closed surface obtained by capping off the puncture.  We may as well assume that $h_-(V) \cap V = \emptyset$, and, following Corollary \ref{cor:isotopic2}, $h(V)$ is not isotopic to $V$ in $F$.  In particular, $V$ must be non-separating.   Since $V$ and $h(V)$ are non-separating and disjoint in $F_-$, but not isotopic in $F$, if $F_-$ is compressed along both $V$ and $h(V)$ simultaneously, $F_-$ becomes a disk.  Apply the Surgery Principle Lemma \ref{lemma:surgprin} to $V$ and conclude that $U \subset S^3$ bounds a disk after $0$-framed surgery on $V$.  In particular, if $N$ is the $3$-manifold obtained by $0$-framed surgery on $V$ alone, then surgery on $U \cup V$ would give $N \# S^1 \times S^2$.  For this to be $\#_{2} (S^{1} \times S^{2})$ would require $N \cong S^1 \times S^2$ hence, again by classical Property R, $V \subset S^3$ would be the unknot.  
\end{proof}

Return to the general case of fibered manifolds and surgery on a curve $c$ in the fiber, and consider the case in which the fiber has genus two.  According to Corollary \ref{cor:isotopic2}, if the result of surgery on $c$ is $\#_{2} (S^{1} \times S^{2})$, then $h(c)$ is not isotopic to $c$ in $F$.   The following Proposition is a sort of weak converse.  

\begin{prop}  \label{prop:nonisotopic} For $M, F, h, c, M_{surg}$ as above, suppose $F$ has genus $2$ and $h(c)$ can be isotoped off of $c$ in $F$.  
If $h(c)$ is not isotopic to $c$ in $F$ then $M_{surg} \cong L \# S^1 \times S^2$, where $L$ is $S^3$, $S^1 \times S^2$, or a Lens space.  
\end{prop}

\begin{proof}  We may as well assume that $h(c)$ is disjoint from $c$ but not isotopic to $c$ in $F$.  Since $F$ is genus two, this  immediately implies that $c$ is non-separating.  

Take the Heegaard viewpoint of Example \ref{example:double}. The complement $M_-$ of a regular neighborhood of  $F'$ in $M_{surg}$ has a Heegaard splitting $H_1 \cup_{F_0} H_2$, with the splitting surface $F_0$ a fiber not containing $c$.  Since $h(c)$ can be isotoped off of $c$ in $F_0$, the Heegaard splitting is a weakly reducible splitting, with $c \subset F_0 = \bdd_+ H_2$ bounding a disk in $H_2$ and $h(c) \subset F_0$ bounding a disk in $H_1$.  

Now do a weak reduction of this splitting.  That is, consider the $2$-handles $C_2 \subset H_2$ with boundary $c \subset F_0 = \bdd_+ H_2$ and $C_1 \subset H_1$ with boundary $h(c)$ in $F_0$ in $N$. Since $c$ and $h(c)$ are disjoint, $N$ can also be regarded as the union of compression bodies $H'_2 = H_2 - C_2 \cup C_1$ and $H'_1 = H_1 - C_1 \cup H_2$.  Each $H'_i$ can be regarded as obtained from $F' \times I$ by attaching a single $2$-handle.  Moreover it is $\bdd_- H'_1$ that is identified with $\bdd_- H'_2$ to get $M_-$.   A genus count shows that this new surface $F'' = \bdd_- H'_i \subset M_-$ is a sphere.  Put another way, the manifold $L = M_{surg} - \eta(F'')$ is Heegaard split by the torus $F'$, so $M_{surg}$ is the connected sum of $S^1 \times S^2$ with a manifold $L$ that has a genus one Heegaard splitting.  \end{proof}

\begin{cor}  \label{cor:nonisotopic} For $M, F, h, c, M_{surg}$ as above, suppose $F$ has genus $2$ and $M_{surg}$  is reducible.  

If $h(c)$ is not isotopic to $c$ in $F$ then $M_{surg} \cong L \# M'$, where $L$ is $S^3$, $S^1 \times S^2$, or a Lens space and $M'$ is either $S^1 \times S^2$ or a torus bundle over the circle.  
\end{cor}
\begin{proof}  Via \cite[Corollary 4.2]{ST} and the proof of Theorem \ref{thm:main} we know that the surgery on $c$ must use the framing given by the fiber in which it lies.  Apply Proposition \ref{prop:monodromy}.  If the first conclusion holds, and $h(c)$ can be isotoped off of $c$ in $F$, then  Proposition \ref{prop:nonisotopic} can be applied and that suffices.  If the second conclusion holds then $c$ is non-separating, so $F'$ is a torus, as required. \end{proof}

\begin{cor}  \label{cor:genustwo} Suppose $U \subset S^3$ is a genus two fibered knot and $V \subset S^3$ is a disjoint knot.  Then $0$-framed surgery on $U \cup V$ gives $\#_{2} (S^{1} \times S^{2})$  if and only if  after possible handle-slides of $V$ over $U$, 
\begin{enumerate}
\item $V$ lies in a fiber of $U$;
\item in the closed fiber $F$ of the manifold $M$ obtained by $0$-framed surgery on $U$, $h(V)$ can be isotoped to be disjoint from $V$;
\item $h(V)$ is not isotopic to $V$ in $F$; and
\item the framing of $V$ given by $F$ is the $0$-framing of $V$ in $S^3$.
\end{enumerate}.  
\end{cor}

\begin{proof}  Suppose first that $0$-framed surgery on $U \cup V$ gives $\#_{2} (S^{1} \times S^{2})$.  Apply  \cite[Corollary 4.2]{ST} as in the proof of Theorem \ref{thm:main} to handle-slide $V$ over $U$ until it lies in the fiber of $U$ in a way that the $0$-framing on $V$ is the framing given by the fiber in which it lies.  Proposition \ref{prop:monodromy} shows that $h(V)$ satisfies the second condition and Corollary \ref{cor:isotopic2} gives the third: $h(V)$ is not isotopic in $F$ to $V$.

For the other direction, suppose $V$ lies in a fiber of $U$ and the four conditions are satisfied.  The last condition says that the surgery on $V$ is via the slope of the fiber.  By Proposition \ref{prop:nonisotopic}, the surgery gives $L \# S^1 \times S^2$, for $L$ either $S^3$, a Lens space, or $S^1 \times S^2$.  But $U$ and $V$ are unlinked in $S^3$ (push $V$ off of $F$), so $0$-framed surgery on $U \cup V$ must give a homology $\#_{2} (S^{1} \times S^{2})$.  This forces $L$ to be a homology $S^1 \times S^2$, hence $S^1 \times S^2$ precisely. \end{proof}

\section{Connected sums of fibered knots}  \label{sect:sums}

There is a potentially useful addendum to Proposition \ref{prop:monodromy} in the case that $h$ has a separating curve that is invariant under the monodromy.  Suppose, as usual, $M$ is an orientable closed $3$-manifold that fibers over the circle, with fiber $F$ and monodromy $h: F \to F$.  Suppose further that there is a separating simple closed curve $\gamma \subset F$, with complementary components $F_1, F_2$ in $F$, so that $h(\gamma) = \gamma$ and $h(F_i) = F_i, i = 1, 2$. Let $h_i = h|F_i: F_i \to F_i$.  

\begin{prop}  \label{prop:arcmono}  Suppose $c \subset F$ is a simple closed curve so that $0$-framed surgery on $c$ in $M$ creates $\#_{2} (S^{1} \times S^{2})$ and $c$ has been isotoped in $F$ to minimize $|c \cap \gamma|$.   For any element $[\ddd] \in H_1(F_i)$ represented by an arc component $\ddd$ of $c - \gamma$, the algebraic intersection satisfies $$-1 \leq [\ddd] \cdot h_{i_*}([\ddd]) \leq 1.$$ \end{prop}

\begin{proof}  Recall the following standard fact about curves in surfaces:  Suppose $\aaa, \beta, \gamma$ are simple closed curves in a surface $F$ so that neither $|\gamma \cap \aaa|$ or $|\gamma \cap \bbb|$ can be reduced by an isotopy of $\aaa$ or $\bbb$.  Then there is an isotopy of $\aaa$ in $F$ that minimizes $|\aaa \cap \beta|$ via an isotopy that never changes $|\gamma \cap \aaa|$ or $|\gamma \cap \bbb|$.  

Apply this fact to the curves $\gamma, \aaa = c, \bbb = h(c)$ in $F$.  Since $M_{surg} \cong \#_{2} (S^{1} \times S^{2})$, the second conclusion of Proposition \ref{prop:monodromy} does not hold, so the first conclusion must:  $h(c)$ can be isotoped to be disjoint from $c$ in $F$.  It follows that for any arc component $\ddd \subset F_i$ of $c - \gamma$, there is a proper isotopy of $h_i(\ddd)$ in $F_i$ so that after the isotopy $\ddd$ and $h(\ddd)$ are disjoint.  Connect the ends of $\ddd$ in $\gamma$ by a subarc of $\gamma$ to get a simple closed curve representing  $[\ddd] \in H_1(F_i)$.  Similarly connect the ends of $h(\ddd)$ to get a representative of $h_{i_*}([\ddd]) \in H_1(F_i)$.  Depending on whether the pairs of ends are interleaved in $\gamma$, the resulting circles can be isotoped either to be disjoint or to intersect in a single point.  \end{proof}

Proposition \ref{prop:arcmono} might give useful information about the monodromy of a connected sum of fibered knots.  Suppose the knot $K \subset S^3$ is the connected sum of two fibered knots $K_1$ and $K_2$.  Then the knot $K$ is also fibered; its fiber is the $\bdd$-connected sum of the fibers for $K_1$ and $K_2$.   This structure carries over to the manifold $M$ obtained by $0$-framed surgery on $K$.  If each $M_i, i = 1, 2$ is the manifold obtained by $0$-framed surgery on $K_i$, with monodromy $h_i: F_i \to F_i$, then the fiber of $M$ is $F = F_1 \# F_2$ and the monodromy $h: F \to F$ is the connected sum of the $h_i$ along an invariant simple closed curve $\gamma \subset F$.  The closed curve $\gamma$ intersects the fiber $F_- \subset S^3$ of $K$ in an invariant arc $\gamma_-$.  The arc $\gamma_-$ can be viewed as the intersection of a decomposing sphere of $K_1 \# K_2$ with the fiber $F_-$.

With this as background, Proposition \ref{prop:arcmono} immediately yields:

\begin{cor} \label{cor:sum}  Suppose $K = K_1 \# K_2$ is a fibered composite knot, and $0$-framed surgery on a link $K \cup V$ creates $\#_{2} (S^{1} \times S^{2})$.  Following \cite{ST}, handle-slide $V$ over $K$ and isotope $V$ so that $V$ lies in a fiber $F_- \subset S^3$ of $K$ and, among all such positionings, choose one that minimizes $|V \cap \gamma_-|$ in $F$.   For  any element $\ddd \in H_1(F_i)$ represented by a component of $V - \gamma_-$, $$-1 \leq [\ddd] \cdot h_{i_*}([\ddd]) \leq 1.$$
\end{cor}

If the summand $K_i$ is a genus one knot, this puts severe restrictions on the set of possible curves of $V \cap F_i$.  For example, suppose $K_i$ is a figure-eight knot.   Then there is a symplectic basis on $H_1(F_i)$ with respect to which  the monodromy $h_{i_*}$ is given by the matrix \[\left( \begin{array}{cc} 2 & 1 \\ 1 & 1 \end{array}\right) : \mathbb{Z}^2 \to \mathbb{Z}^2.\]
For a class $\left( \begin{array}{c} m \\ n \end{array}\right) \in H_1(F_i) \cong \mathbb{Z}^2 $ to have the property 
\[-1 \leq \left( \begin{array}{cc} m & n \end{array} \right) 
 \left( \begin{array}{cc} 0 & -1 \\ 1 & 0 \end{array} \right)
 \left( \begin{array}{cc} 2 & 1 \\ 1 & 1 \end{array} \right) \left( \begin{array}{c} m \\ n \end{array}\right) \leq 1\]
requires \[ -1 \leq -m^2 + mn + n^2 \leq 1.\]  An elementary descent argument shows that solutions are pairs $(m, n)$ such that $m \cdot n \geq 0$ and $|n|, |m|$ are successive Fibonacci numbers or $m\cdot n \leq 0$ and $|m|, |n|$ are successive Fibonacci numbers.  As many as three of these curves may be present simultaneously: if $f_1, f_2, f_3, f_4$ are successive Fibonacci numbers, then a similar calculation shows that the three successive pairs \[\left( \begin{array}{c}  f_1 \\ f_2 \end{array}\right),  
\left( \begin{array}{c} f_2 \\  f_3 \end{array}\right) , 
  \left( \begin{array}{c}  f_3 \\  f_4 \end{array}\right) \]
 in $\mathbb{Z}_2 \cong H_1(F_i)$ may be represented in the punctured torus by disjoint arcs. 

Similarly, for the trefoil knot, there is a symplectic  basis on $H_1(F_i)$ with respect to which the monodromy $h_{i_*}$ is given by the matrix \[\left( \begin{array}{cc} 0 & 1 \\ -1 & 1 \end{array}\right)\]
For a class $\left( \begin{array}{c} m \\ n \end{array}\right) \in H_1(F_i) \cong \mathbb{Z}^2 $ to have the property \[-1 \leq \left( \begin{array}{cc} m & n \end{array} \right)
 \left( \begin{array}{cc} 0 & -1 \\ 1 & 0 \end{array} \right)
  \left( \begin{array}{cc} 0 & 1 \\ -1 & 1 \end{array} \right) \left( \begin{array}{c} m \\ n \end{array}\right) \leq 1\]
requires \[ -1 \leq m^2 + mn + n^2 \leq 1. \]  This allows only three possible curves: $$  \left( \begin{array}{c} m \\ n \end{array}\right) = 
 \left( \begin{array}{c} \pm 1 \\ 0 \end{array}\right),  \left( \begin{array}{c} 0 \\ \pm 1 \end{array}\right) , {\rm or}  \left( \begin{array}{c} \pm 1 \\ \mp 1 \end{array}\right) .$$
 
These three can be represented in the punctured torus by disjoint arcs.

\bigskip

{\bf Afterword:}  It will be shown in \cite{GSch} that a similar analysis gives a precise catalog of all possible curves $V$ in the complement of a square knot $K$ (up to band-sum with $K$) so that surgery on $K \cup V$ gives $\#_{2} (S^{1} \times S^{2})$.  It will also be shown that the central example of \cite{Go} gives rise to a likely counterexample to Property 2R in which one of the link components is the square knot.  The other component can be precisely described, but it remains a puzzle how, even after band-sums with the square knot, it fits into this catalog.  So it also remains mysterious how, via \ref{thm:main} it then gives rise to a probable genus one non-fibered counterexample to Property 2R.


\begin{thebibliography}{5}

\bibitem[CG]{CG} A.~Casson and C.~McA.~Gordon, Reducing Heegaard splittings,  {\em Topology  and its applications}, {\bf 27} (1987) 275--283.

\bibitem[Fr]{Fr}  M.~Freedman, The topology of four-dimensional manifolds.   {\em J. Differential Geom.  }, {\bf 17}  (1982) 357--453. 

\bibitem[Ga1]{Ga1} D.~Gabai, Foliations and the topology of $3$-manifolds. II. {\em J. Differential Geom.} {\bf 26} (1987) 461--478.

\bibitem[Ga2]{Ga2} D.~Gabai, Foliations and the topology of $3$-manifolds. III, {\em J. Differential Geom.} {\bf 26} (1987) 479--536.

\bibitem[Ga3]{Ga3} D.~Gabai,  Surgery on knots in solid tori {\em Topology} {\bf 28} (1989) 1-6.

\bibitem[GSch]{GSch} R.~Gompf and M.~Scharlemann, Does the square knot have Property 2R? {\it to appear}.

\bibitem[GS]{GS} R.~Gompf and A.~Stipsicz, \underline{$4$-Manifolds and Kirby Calculus}, Graduate Studies in Mathematics {\bf 20} American Mathematical Society, Providence, RI, 1999.

\bibitem[Go]{Go} R.~Gompf, Killing the Akbulut-Kirby $4$-sphere, with relevance to the Andrews-Curtis and Schoenflies problems {\em Topology}  {\bf 30}  (1991),  97--115.

\bibitem[Ki1]{Ki1} R Kirby, A calculus for framed links in $S^3$ {\em Invent. Math} {\bf 45} (1978) 35--56. 

\bibitem[Ki2]{Ki2} R Kirby, , Problems in Low-Dimensional Topology, in {\underline Geometric Topology}, Edited by H. Kazez, AMS/IP Vol. 2, International Press, 1997.

\bibitem[ST]{ST}  M.~Scharlemann, A.~Thompson, Surgery on knots in surface $\times$ I, {\it ArXiv 0807.0405}

\end{thebibliography}
\end{document}